\documentclass[a4paper,12pt]{amsart}

\usepackage{amsmath,amsthm,amsfonts,amssymb,url}
\usepackage[T1]{fontenc}
\usepackage[latin1]{inputenc}
\usepackage[francais]{babel}
\usepackage{xypic}

\newcommand{\jj}{\mathrm{\got j}}
\newcommand{\CC}{\mathcal{C}}
\newcommand{\ssi}{si, et seulement si, }
\newcommand{\ord}{\mathrm{ord}}
\newcommand{\eis}{\mathrm{eis}}
\newcommand{\ps}{\par \smallskip }
\newcommand{\pn}{\par \noindent}
\newcommand{\Q}{\mathbb{Q}}
\newcommand{\Z}{\mathbb{Z}}

\newcommand{\W}{\mathcal{W}}
\newcommand{\OO}{\mathcal{O}}
\newcommand{\HH}{\mathcal{H}}
\newcommand{\End}{\mathrm{End}}
\newcommand{\Spec}{\mathrm{Spec}}
\newcommand{\Frob}{\rm Frob}
\newcommand{\G}{\mathbb G}
\newcommand{\GL}{\mathrm{GL}}
\newcommand{\Hom}{\mathrm{Hom}}
\newcommand{\Qpb}{\overline{\Q}_p}

\newcommand{\Gal}{\mathrm{Gal}}

\newcommand{\N}{\mathbb{N}}
\newcommand{\pf}{ {\it Preuve: } }
\newcommand{\epf}{$\square$ \par }
\newcommand{\rig}{\mathrm{rig}}
\newcommand{\A}{\mathbb{A}}
\newcommand{\Ext}{\mathrm{Ext}}

\newcommand{\Ker}{\rm Ker}
\newcommand{\tr}{\mathrm{tr}}
\newcommand{\Max}{\mathrm{Specmax}}

\newtheorem{thm}{Th\'eor\`eme}
\newtheorem{prop}{Proposition}
\newtheorem{cor}{Corollaire}
\newtheorem{lemme}{Lemme}

\def\]{\textup{\mbox{]\hspace{-.15em}]}}}
\def\[{\textup{\mbox{[\hspace{-.15em}[}}}

\def\got{\mathfrak}

\title{Lissit\'e de la courbe de Hecke de $\GL_2$ aux points Eisenstein
critiques.}

\author{J. Bellaïche et G. Chenevier}

\begin{document}

\maketitle

\medskip

{\it Abstract}: Let $p$ be a prime number and $\CC$ be the $p$-adic tame
level $1$ eigencurve introduced by Coleman-Mazur. We prove that $\CC$ is
smooth at the evil Eisenstein points and we give necessary and sufficient
conditions for etaleness of the map to the weight space at these points in terms of
$p$-adic zeta values. A key step is the determination at these points of the schematic
reducibility locus of the pseudo-character carried by $\CC$
restricted to a decomposition group at $p$. Then, the smoothness appears to
be a consequence of the fact that the Dirichlet $L$-functions only have simple zeros
at integers. \ps

\medskip

\section*{Introduction}

       Soient $p$ un nombre premier et $\CC$ la courbe de Hecke $p$-adique
de $\GL_2$ de niveau mod\'er\'e $1$, "{\it the eigencurve}", introduite par
Coleman et Mazur dans \cite{eigen}. Notons $G$ le groupe de Galois de $\Q$ non ramifi\'e hors de
$p$ et $T: G \rightarrow A(\CC)$ le pseudo-caract\`ere de dimension $2$
port\'ee par $\CC$. Les points Eisenstein critiques sont les points
$x_{(k,\varepsilon)}$ de $\CC$ en lesquels $U_p(x_{(k,\varepsilon)})=p^{k-1}$
et $T$ est la trace de la repr\'esentation
$$\Q_p.\varepsilon \oplus \Q_p(1-k),$$
$k\geq 2$ \'etant un entier et $\varepsilon$ un caract\`ere d'ordre fini de
$G$ non trivial si $k=2$ et v\'erifiant $\varepsilon(-1)=(-1)^k$. 
Ces points forment une partie discr\`ete de la
partie non ordinaire de $\CC$ et tout point de cette derni\`ere o\`u $T$ est
r\'eductible est de cette forme. Dans cet article, nous prouvons que $\CC$ est lisse en
$x_{(k,\varepsilon)}$ et nous donnons un crit\`ere pour que le morphisme
vers l'espace des poids, $\kappa: \CC \rightarrow \W$, y soit \'etale. \ps

On  conjecture que $\CC$ est lisse en tous ses points classiques.
Aux points classiques r\'eductibles\footnote{Nous dirons qu'un point $x$ de
$\CC$ est r\'eductible (resp. irr\'eductible) si l'\'evaluation de $T$ en $x$ l'est.} ordinaires, il est ais\'e de v\'erifier
que $\CC$ est lisse et m\^eme que $\kappa$ est \'etale, de sorte que nous d\'emontrons
cette conjecture pour tout point classique r\'eductible. En ce qui concerne la
lissit\'e aux points classiques irr\'eductibles, elle est connue dans de
nombreux cas par les travaux de Kisin \cite[thm. 11.10]{kis}.\ps

Notons \'egalement qu'en un point classique non critique $z$, le th\'eor\`eme de
classicit\'e de Coleman montre que le degr\'e de $\kappa$ en $z$ est \'egal
\`a la dimension de l'espace caract\'eristique (pour les $T_l$, $l\neq p$,
et $U_p$) de la forme $f$ associ\'ee
\`a $z$ dans l'espace des formes classiques de m\^eme poids que $f$. En
particulier, dans le cas non critique, $\kappa$ est \'etale en $z$ \ssi
$U_p$ agit de mani\`ere semi-simple sur cet espace caract\'eristique, ce qui
est conjectur\'e, et implique la lissit\'e dans certains cas (cf. \cite[cor. 7.6.3]{eigen}).
\ps 

Notre d\'emonstration repose sur l'\'etude des lieux de r\'eductibilit\'e sch\'ematique $\Spec(R)$ et $\Spec(R_p)$ respectifs de
$T$ et $T_{|D}$ au voisinage d'un point $x$ Eisenstein
critique, $D \subset G$ \'etant un groupe de d\'ecomposition en $p$. Nous montrons en utilisant un r\'esultat de Kisin que $\Spec(R_p)$
est inclus dans la fibre en $x$ de $\kappa$, puis qu'ils sont \'egaux en
utilisant le cas limite du crit\`ere de classicit\'e de Coleman. Nous en
d\'eduisons que $\Spec(R)$ est le point ferm\'e r\'eduit $x$. En  
utilisant des techniques de Mazur-Wiles \cite{MW}, nous obtenons une
majoration du nombre minimal de g\'enerateurs de l'id\'eal d\'efinissant $\Spec(R)$ (donc ici de l'id\'eal maximal $m$ de $\OO_{\CC,x}^{\rig}$) en terme de la dimension de certains
groupes de Selmer. La principalit\'e de $m$ appara\^it alors comme
cons\'equence de ce que les fonctions $L$ de Dirichlet n'ont que des
z\'eros simples aux entiers et des conjectures de Bloch-Kato, connues pour ces derni\`eres. En ce qui concerne 
le degr\'e de $\kappa$ en $x$, nombre assez myst\'erieux du fait que $x$ est
critique, nous montrons qu'il vaut $1$ \ssi une certaine valeur explicite de
la fonction $\zeta$ $p$-adique est non nulle. Cette non-annulation est   
conjectur\'ee, mais non connue; elle est \'equivalente \`a la non-annulation
d'un certain r\'egulateur $p$-adique.  \ps

\medskip

{\small Les auteurs sont heureux de remercier Pierre Colmez,
Barry Mazur, Christophe Soul\'e et Jacques Tilouine pour leur soutien et des
discussions utiles, ainsi que le C.I.R.M.
(Luminy, France) o\`u une partie de ce travail
a \'et\'e r\'ealis\'ee en novembre 2003. L'un des auteurs (J. B.) 
remercie de plus l'I.P.D.E. pour son soutien financier et l'université de 
Rome I pour son hospitalité.} \ps

\medskip

{\small {\sc Notations}}: $p$ est un nombre premier, $\overline{\Q}_p$ une
cl\^oture alg\'ebrique fix\'ee de $\Q_p$, et $v: \Qpb^* \rightarrow \Q$ la
valuation $p$-adique normalis\'ee par $v(p)=1$. Si $X/\Q_p$ est un espace rigide, 
on note $A(X)$ l'anneau des fonctions analytiques globales sur $X$ et
$\OO_X^{\rig}$ le faisceau structural de $X$. Nous 
entendrons par $X(\Qpb)$ la r\'eunion des $X(F)$ o\`u $F$ parcourt les sous-extensions
finies de $\Qpb$.

\section{Rappels sur $\CC$}

La r\'ef\'erence pour cette partie est \cite{eigen}.\ps

\subsection{} Soient $p$ un nombre premier,
$\W:=\Hom_{{\rm gr-cont}}(\Z_p^*,\G_m^{rig})_{pairs}$, et $\kappa: \CC \rightarrow \W$
la courbe de Hecke $p$-adique de niveau mod\'er\'e $1$ pour le groupe $\GL_2$. 
C'est la courbe analytique\footnote{Pr\'ecis\'ement, la
courbe $\CC$ \'etudi\'ee ici est celle not\'ee $D$ loc.cit. Nous
n'utiliserons pas l'identification de $D$ avec la nilr\'eduction de l'espace $C_p$
consid\'er\'e aussi loc.cit.} sur $\Q_p$ construite par Coleman et Mazur
dans \cite{BMF}  et \cite[Chap. 7]{eigen} (voir aussi \cite{B} pour $p=2$) \`a partir du syst\`eme de modules de Banach
sur $\W$ des formes modulaires $p$-adiques surconvergentes. 
La courbe $\CC$ est
s\'epar\'ee, r\'eduite, \'equidimensionnelle de dimension $1$. Le morphisme $\kappa$ est
plat, localement fini. \ps
 
\subsection{}\label{preparation} Soit $\HH:=\Z[\{T_l, l\neq p\},U_p]$. On dispose par construction (\cite[Chap. 7]{eigen}) d'un morphisme d'anneaux 
$\HH \rightarrow A(\CC)$  de sorte que
l'on verra les \'el\'ements de $\HH$ comme des fonctions analytiques
globales, born\'ees par $1$ partout, sur $\CC$.
Par construction toujours (\cite{eigen}), l'application canonique "syst\`eme de valeurs propres"
$\chi: \CC(\Qpb) \rightarrow \mathrm{Hom}_{\rm ann}(\HH,\Qpb)$ est injective, et
identifie $\CC(\Qpb)$ \`a l'ensemble des formes modulaires $p$-adiques
surconvergentes propres, de pente finie, et de niveau mod\'er\'e $1$. Soient $F/\Q_p$ un corps local, $x \in \CC(F)$. On 
note $f_x$ l'unique forme $p$-adique surconvergente propre normalis\'ee
correspondante\footnote{Cela existe toujours sauf si $x$ est sur la droite
Eisenstein ordinaire et de poids trivial, auquel cas on pose
$f_x=1$ (cf. \cite[prop. 3.6.1]{eigen}).}, et $M(x) \subset M_{\kappa(x)}^{\dagger}$ l'espace
caract\'eristique pour $\HH$ de $f_x$. L'image $\HH(x)$ de $F \otimes_{\Z} \HH$ dans
$\End_F(M(x))$ est une $F$-alg\`ebre locale de dimension finie. 
L'accouplement standard $M(x) \times \HH(x) \rightarrow F$, $(f,h) \mapsto
a_1(h(f))$ est non d\'eg\'en\'er\'e, sauf si $\kappa(x)=1$ et $x$ est sur la
droite Eisenstein ordinaire, auquel cas il est nul mais 
$\dim_F(M(x))=\dim_F(\HH(x))=1$ (cf.
\cite[prop. 3.6.1]{eigen}). \ps

\begin{prop} \label{voisinage} Si $x \in \CC(F)$, il existe un $F$-voisinage affinoide $\Omega$ de $x$
tel que: \ps

(a) $\kappa(\Omega)$ est un ouvert affinoide,\ps
(b) $\kappa_{|\Omega}$ est fini et plat, \'etale hors de $x$, de degr\'e
$\dim_F M(x)$, \ps
(c) la fibre de $\kappa_{|\Omega}$ au dessus de $\kappa(x)$ s'identifie
canoniquement \`a $\Spec(\HH(x))$. \pn
De plus, l'application naturelle $\OO_{\W,\kappa(x)}^{\rig} \otimes_{\Z}\HH \rightarrow
\OO_{\CC,x}^{\rig}$ est surjective.\ps

\end{prop}

\pf Choisir $\Omega$ tel que (a) et le premier point de (b) soient vrai est
possible par construction. Posons $V:=\kappa(\Omega)$, $F_x$ la fibre de
$\kappa_{|\Omega}$ au dessus de $\kappa(x)$. D'apr\`es \cite[lemme 2.1.6]{Be}, les
$\kappa^{-1}(U)$ avec $\kappa(x) \in U$ forment une base de voisinages rigides
analytiques de $F_x$. Ainsi, quitte \`a r\'eduire $V$ et remplacer $\Omega$
par sa composante connexe contenant $x$, on peut supposer que $F_x$ est un
sch\'ema local, et que $\kappa$ est \'etale hors de $x$. Soit $y \in V
\backslash \{\kappa(x)\}$, alors par construction, le degr\'e de $\kappa$ est
$|\kappa^{-1}(y)|=\sum_{\kappa(z)=y} \dim_F(\HH(y))=\sum_{\kappa(z)=y}\dim_F(M(y))$. Ce 
degr\'e est aussi $\dim_F(M(x))$ car la famille de formes modulaires d\'ecoup\'ee par 
$\Omega$ est localement libre sur $A(V)$, ce qui termine de prouver (b). 
En r\'eappliquant le raisonnement ci-dessus au point $x$
et en utilisant la platitude de $\kappa$ et $\dim_F(M(x))=\dim_F(\HH(x))$,
on en d\'eduit (d). La derni\`ere assertion est satisfaite par construction. \epf


\subsection{} \label{composantes}Terminons cette section par la description de certaines
composantes particuli\`eres de $\CC$. Le {\it lieu parabolique} de $\CC$,
que l'on notera $\CC^0$, est le ferm\'e r\'eduit de $\CC$ d\'efini par 
$$\CC^0:=\{x \in \CC, \, \, f_x \textrm{ s'annule \`a la pointe }
\infty\}.$$ 
Par construction de $\CC$, c'est un ferm\'e Zariski de $\CC$ qui est
d'\'equidimension $1$. La restriction de $\kappa$ \`a $\CC^0$ est encore fini et
plate. \par
	Enfin, le {\it lieu ordinaire} de $\CC$, que l'on notera $\CC^{\rm
ord}$, est l'ouvert admissible de $\CC$ d\'efini par 
$$\CC^{\rm ord}:=\{x \in \CC, \, \, |U_p(x)|=1\}.$$
La relative compacit\'e de l'image de $\HH$ dans $A(\CC)$ (cf. par
exemple \cite[\S 4.6 rem. i.]{JLp}) assure que 
l'idempotent de Hida $e:=\lim_{n \rightarrow \infty} U_p^{n!}$ d\'efinit un \'el\'ement de
$A(\CC)$. Cela montre que $\CC^{\ord}$ est en fait l'ouvert ferm\'e admissible
de $\CC$ d\'efini par $e=1$.
On pourrait d\'emontrer que $\kappa: \CC^{\ord} \rightarrow \W$ est fini, et que
c'est la fibre g\'en\'erique de l'alg\`ebre de Hecke ordinaire de Hida, 
mais nous n'en aurons pas besoin.

\section{Rappels sur la th\'eorie de Mazur-Wiles}

Nous rappelons dans cette section, en les \'etendant l\'eg\`erement, 
quelques r\'esultats d\'emontr\'es dans \cite{MW} (voir aussi
\cite{HP}).

\subsection{}\label{MWintro} Soit $(A,m,k)$ un anneau local noeth\'erien
hens\'elien r\'eduit. On note $K=\prod_j K_j$
son anneau total de fractions, et on suppose donn\'ee $\rho=(\rho_j): G \rightarrow
\GL_2(K)$ une repr\'esentation de trace not\'ee $T$ telle que $T(G) \subset A$, de
d\'eterminant dans $A$. On suppose que $T \bmod m$ est somme de deux
caract\`eres {\it distincts} $\chi_i: G \rightarrow k^*$, $i=1,\, 2$. \ps

Fixons $s \in G$ tel que $\chi_1(s) \neq \chi_2(s)$. Le
polyn\^ome caract\'eristique de $s$ est scind\'e dans $A$ car $A$ est
hens\'elien, \`a racines distinctes dans chacun des $K_j$. On note
$\lambda_i \in A$ l'unique racine telle que $\lambda_i \bmod m= \chi_i(s)$. On
peut donc trouver une $K$-base de $K^2$, $e_1,e_2$ telle que
$s(e_i)=\lambda_i e_i$. Une telle base sera dite {\it adapt\'ee \`a $s$}.
On note $a,b,c,d$ les coefficients matriciels de $\rho$ dans cette base,
et $B$ et $C$ les sous-$A$-modules de $K$ engendr\'es par les $b(g)$ et
$c(g')$ respectivement.

\subsection{} Soit $I \subsetneq A$ un id\'eal tel que $T \bmod I$ soit la
somme de deux caract\`eres $\psi_1, \, \psi_2: G \rightarrow (A/I)^*$, tels que $\psi_i \bmod m = \chi_i $.

\begin{prop} \label{MW} $\Hom_A(B,A/I)$ s'injecte dans
$\Ext^1_{(A/I)[G]}(\psi_2,\psi_1)$. De plus, \ps

\noindent (a) si les $\rho_j$ sont semi-simples, $B$ est un sous-$A$-module
de type fini de $K$, \ps

\noindent (b) si les $\rho_j$ sont irr\'eductibles, alors l'annulateur de
$B$ est nul.

\end{prop}

\medskip

\begin{lemme} \label{idempotents} Pour tout $g \in G$, on a $a(g),\, d(g)
\in A$ et $a(g)-\psi_1(g), \, d(g)-\psi_2(g) \in I$. De plus, pour tous $g,\, g' \in
G$, $b(g)c(g') \in I$.
\end{lemme}

\pf Les \'el\'ements $T(sg)=\lambda_1 a(g)+ \lambda_2 d(g)$ et
$T(g)=a(g)+d(g)$ sont dans $A$, ainsi donc que $a(g)$ et $d(g)$ car $\lambda_1-\lambda_2$ est
inversible dans $A$. En r\'eduisant modulo $I$ les deux relations plus haut, il vient que
$a(g)-\psi_1(g)$ et $d(g)-\psi_2(g)$ sont solutions du syst\`eme $x+y=0$ et
$\bar{\lambda}_1x+\bar{\lambda}_2y=0$ qui est inversible dans $A/I$, car
dans $A/m$. Cela conclut le premier point. Le second point en d\'ecoule, car
\begin{equation} \label{produit} a(gg')=a(g)a(g')+b(g)c(g'). \end{equation}
\epf

\smallskip
\medskip

Notons $\bar{b}$ l'image de $b$ dans $B/IB$. Un cons\'equence imm\'ediate du
lemme \ref{idempotents} est le:

\begin{lemme} \label{extension} L'application 
$$ G \longrightarrow \left( \begin{array}{cc} (A/I)^*  &  B/IB 
\\ 0 & (A/I)^* \end{array}\right), \, \, \, \, g \mapsto \left(
\begin{array}{cc} \psi_1(g) &  \overline{b(g)} \\ 0 & \psi_2(g) \end{array}
\right), $$
est un morphisme de groupes. \end{lemme}


En particulier, on dispose d'une une application $A$-lin\'eaire $$\jj:
\Hom_A(B/IB,A/I) \rightarrow \Ext^1_{(A/I)[G]}(\psi_2,\psi_1).$$ \noindent Posons $H:=\ker(\psi_1/\psi_2)$.
\ps

\begin{lemme} $B/IB$ est engendr\'e comme $A$-module par les $b(h), h \in
H$, et $\jj$ est injective.
\end{lemme}

\pf Soit $g$ dans $G$, un calcul montre que $$b(sgs^{-1}g^{-1})=\frac{b(g)}{\psi_2(g)}(\frac{\lambda_1}{\lambda_2}-1)$$
dans $B/IB$. Comme $sgs^{-1}g^{-1} \in H$, cela conclut le premier point. Soit $f \in
\Ker(\jj)$, $g \mapsto f(b(g))$ est un cobord, donc trivial restreint \`a $H$. Donc $f$
est nulle sur $A.b(H)=B/IB$ par le premier point. 
\epf

\begin{lemme} \label{tf} Si les $\rho_j$ sont semi-simples, $B$ est un
$A$-module de type fini. \end{lemme}

\pf Comme le $A$-module $B$ est un quotient de $A[\rho(G)]$, il suffit de
montrer que ce dernier est de type fini sur $A$. Comme $A[\rho(G)] \subset
\prod_j \rho_j(A[G])$, on peut supposer que $K$ est un corps. Comme $\rho$
est semi-simple, la trace de $M_2(K)$ est non d\'eg\'en\'er\'ee sur
$K[\rho(G)]$. De plus, elle est \`a valeurs dans $A$ sur le sous-$A$-module
$R:=A[\rho(G)]$. Comme $A$ est noeth\'erien et $R$ engendre $K[\rho(G)]$
comme $K$-espace vectoriel, un argument standard montre que $R$ est de type
fini. \epf
\medskip

Enfin, il est clair que si $\rho_j$ est irr\'eductible, $\mathrm{Im}(B \rightarrow
K_j)$ est non nulle. Cela ach\`eve la preuve de la proposition \ref{MW}. $\square$ \ps

\subsection{} \label{defidred} Remarquons que
le lemme \ref{idempotents} montre que $BC
\subset m$ est le plus grand id\'eal $J$ de $A$ ayant la propri\'et\'e
que $T \bmod J$ est somme de deux caract\`eres. On l'appellera {\it
l'id\'eal de r\'eductibilit\'e de $T$}. \ps

\medskip

\begin{cor} \label{corMW} Supposons les $\rho_j$ semi-simples. \ps
(a) Si $\dim_k(\Ext^1_{k[G]}(\chi_2,\chi_1))=1$, alors $\rho$ est d\'efinie
sur $A$. Si de plus l'annulateur de $B$ est nul, $B$ est libre de rang $1$
sur $A$, et il existe une base adapt\'ee \`a $s$ dans laquelle $B=A$.
\ps

(b) Si l'id\'eal de r\'eductibilit\'e de $T$ est l'id\'eal maximal
de $A$ et si $\dim_k(\Ext^1_{k[G]}(\chi_1,\chi_2))=\dim_k(\Ext^1_{k[G]}(\chi_2,\chi_1))=1$,
alors $A$ est de valuation discr\`ete. \end{cor} 

\pf Prouvons le (a). Si l'on applique la proposition \ref{MW} \`a l'id\'eal
maximal de $A$, il vient que $\dim_k B\otimes_A k \leq
\dim_k(\Ext^1_{k[G]}(\chi_2,\chi_1))=1$. Comme $B$ est de
type fini par le lemme \ref{tf}, $\dim_k B \otimes_A k=1$ et
donc $B$ est un $A$-module monog\`ene par Nakayama, ainsi donc que son image
$B_j$ dans $K_j$. Posons $f_j=1$ si $B_j=0$ et $(f_j)=B_j$ sinon, on a $f:=(f_j) \in
K^*$. Ainsi, quitte \`a
remplacer $e_2$ par $f^{-1} e_2$, on conclut le (a). \ps Si pour $(i_1,i_2)=(1,2)$ et $(2,1)$ on a
$\dim_k(\Ext^1_{k[G]}(\chi_{i_1},\chi_{i_2}))=1$, alors pour les m\^emes raisons que
plus haut, $B$ et $C$
sont monog\`enes, ainsi donc que
$m=BC$ par hypoth\`ese. L'anneau $A$ \'etant r\'eduit, il est donc de valuation
discr\`ete. \epf

\subsection{}\label{topo} En vue d'appliquer les r\'esultats de cette
section \`a un groupe topologique, nous avons besoin d'un sorite de
topologie. On conserve les hypoth\`eses du \S \ref{MWintro} et on suppose que $G$ est
un groupe topologique. On suppose de plus que $A$
est un anneau topologique s\'epar\'e ayant la propri\'et\'e suivante: (TOP)
le foncteur
d'oubli des
$A$-modules topologiques s\'epar\'es de type fini vers les $A$-modules de
type fini admet une section pleinement fid\`ele munissant $A$ de sa topologie. On
fixe une telle section, de sorte que tout $A$-module de type fini est muni de la topologie
donn\'ee par cette section. En particulier, si $I$ est un id\'eal de $A$, 
on dispose d'une topologie sur $A/I$. Noter qu'un sous-$A$-module $N$ d'un
$A$-module $M$ de type fini est automatiquement ferm\'e, car $M/N$ est
s\'epar\'e et $M \rightarrow M/N$ continue. \ps

\begin{prop} \label{topopo} Supposons que $T: G \rightarrow A$ est continu et que les
$\rho_j$ sont semi-simples. Alors les $\psi_i: G \rightarrow (A/I)^*$ sont continus et l'application
$$\jj: \Hom_A(B/IB,A/I) \rightarrow
\Ext^1_{(A/I)[G]}(\psi_2,\psi_1)$$ a son
image dans $\Ext^1_{cont, (A/I)[G]}(\psi_2,\psi_1)$. \end{prop}

\pf D'apr\`es le lemme \ref{idempotents} et sa preuve, 
$\psi_1$ co\"incide avec $a \bmod I$ et
$$a(g)=\frac{T(sg)-\lambda_2T(g)}{\lambda_1-\lambda_2}.$$ Ainsi, $g \mapsto
a(g), \psi_1(g)$ et $\psi_1(g)^{-1}=\psi_1(g^{-1})$ sont continus, car
$T$ l'est et par (TOP). Il en va de m\^eme pour $\psi_2$. 
Cela a donc un sens de parler d'extensions continues entre $\psi_2$ et $\psi_1$. 
Comme $B$ est de type fini par le lemme \ref{tf}, et que toute application $A$-lin\'eaire $B \rightarrow A/I$ est
continue par (TOP), il ne reste qu'\`a montrer que $b: G \rightarrow
B$ est continue.  Soit $B_j$ l'image de $B$ dans $K_j$; c'est un
quotient de $B$, donc de type fini sur $A$. L'application canonique $B \rightarrow \prod_j B_j$
est injective et c'est un hom\'eomorphisme sur son image. Il suffit donc de v\'erifier que $g \mapsto b(g)_j$ est
continue. On peut supposer que $K$ est un corps, puis que $\rho$ est
irr\'eductible, car sinon $C=B=0$. Soit $g' \in G$ tel que $c(g') \neq 0$, la multiplication
par $c(g')$ induit un hom\'eomorphisme de $B$ sur son image dans $A$. Il suffit donc de
v\'erifier que $g \mapsto c(g')b(g) \in A$ est continue, ce qui d\'ecoule de
ce que $a$ l'est et de la formule (\ref{produit}). \epf

\medskip
Une preuve identique \`a celle du corollaire \ref{corMW} d\'emontre alors le:

\begin{cor} \label{corMWtop} Sous les hypoth\`eses de la proposition \ref{topopo}, le
corollaire \ref{corMW} reste vrai s'il on remplace dans son \'enonc\'e les
groupes d'extensions mis en jeu par leurs sous-groupes d'extensions continues.
\end{cor}

\medskip

{\it Exemple\label{exemple}:} Soient $k$ un corps local non archim\'edien, $X$ un
$k$-affinoide r\'eduit, $x \in X$ et $A$ l'anneau local rigide en $x$. On rappelle que
l'anneau $A$ est limite inductive filtrante des $A(U)$ o\`u $U$ est un ouvert affinoide de $X$ contenant $x$.
C'est un anneau local noeth\'erien r\'eduit (\cite[\S 7.3.2]{BGR}) et hens\'elien
(\cite[\S 2.1]{Be}). On le munit de la topologie localement convexe la plus fine telle que les
$A(U) \rightarrow A$ soient continues (cf. \cite[ch. I, E]{found}),
$A(U)$ \'etant muni de sa topologie de $k$-espace de Banach. En particulier, si $m$ est l'id\'eal maximal de $A$, la projection canonique $A \rightarrow A/m^n$ est continue.
Cette topologie fait de $A$ une $k$-alg\`ebre topologique, elle est s\'epar\'ee
car $A/m^n$ l'est et $\bigcap_{n\geq 0}m^n=\{0\}$. Si $M$ est un $A$-module de type fini et $f: A^n \rightarrow M$ une
surjection $A$-lin\'eaire, la topologie localement convexe quotient de $M$
le munit d'une structure de $A$-module topologique qui est en
fait ind\'ependante de la surjection $f$ choisie, et s\'epar\'ee. On voit facilement que toute 
application lin\'eaire entre deux $A$-modules de type fini est continue. 
Ainsi, $A$ satisfait (TOP). \ps

\section{La pseudo-repr\'esentation port\'ee par $\CC$.}

\subsection{} Soient $G$ le groupe de Galois de la sous-extension maximale de
$\overline{\Q}$ non
ramifi\'ee hors de $p$, $D \subset G$ un groupe de d\'ecomposition en $p$ attach\'e
\`a un plongement $\overline\Q \hookrightarrow \Qpb$ que l'on fixe. 
On note $Z \subset \CC(\Qpb)$ l'ensemble des points classiques de $\CC$.
Par d\'efinition $z \in Z$ \ssi $\chi(z)$ est le syst\`eme de valeurs
propres d'une forme modulaire classique sur $X_1(p^n)$ pour un certain
entier $n\geq 0$. Les formes modulaires apparaissant ainsi sur $\CC$ sont exactement celles qui ne
sont pas supercuspidales en $p$. On sait que l'ensemble des points
ferm\'es de $\CC$ ainsi obtenu est tr\`es Zariski-dense dans $\CC$. \ps

\subsection{}\label{Tcont} \`A chaque $z \in Z$ est associ\'ee, par les travaux de
Eichler-Shimura, Igusa, Deligne, une
unique repr\'esentation semi-simple continue $$\rho_z: G \rightarrow
\GL_2(\Qpb),$$
ayant la propri\'et\'e que la trace d'un Frobenius g\'eom\'etrique en $l \neq p$ est
$T_l(x)$. La compacit\'e relative de l'image de $\HH$ dans $A(\CC)$, et le fait
que $A(\CC)$ est r\'eduit, entra\^inent
que la trace $x \mapsto \tr(\rho_x)$ de ces repr\'esentations se prolonge analytiquement en
une unique pseudo-repr\'esentation continue de dimension $2$: 
	$$T: G \rightarrow A(\CC),$$ 

\noindent satisfaisant $T(F_l)=T_l$. En particulier, pour tout $x \in \CC(\Qpb)$, il
existe\footnote{En fait, si $x \in \CC(F)$, alors $\rho_x$ est d\'efinie sur
$F$. Cela vient de ce que $\rho_x(\Frob_{\infty})$ a pour polyn\^ome
caract\'eristique $(X-1)(X+1)$, donc l'obstruction \`a ne pas \^etre
d\'efinie sur $F$ est nulle.} une unique repr\'esentation {\it semi-simple} continue $\rho_x: G
\rightarrow \GL_2(\Qpb)$ dont la trace est l'\'evaluation en $x$ de $T$. En
particulier, un point $x \in \CC(\Qpb)$ est uniquement d\'etermin\'e par
le couple $(\rho_x,U_p(x))$. \ps

\subsection{} Si $x \in \CC(\Qpb)$, on sait que le polyn\^ome de Sen de
$(\rho_x)_{|D}$ est $T(T-d\kappa(x)+1)$, o\`u $d\kappa(x)$ d\'esigne la
d\'eriv\'ee en $1$ du caract\`ere de $\Z_p^*$ associ\'e \`a $\kappa(x)$.
Par les travaux de Kisin (\cite[thm. 6.3]{kis}), on a 
$$D_{\mathrm{cris}}((\rho_x)_{|D})^{\varphi=U_p(x)} \neq 0.$$

\section{Points r\'eductibles de $\CC$}

\subsection{Points Eisenstein} Commen\c{c}ons par d\'efinir les points
{\it Eisenstein critiques} de $\CC$. Soient $k\geq 2$ un entier et $\varepsilon: \Z_p^*
\rightarrow \Qpb^*$ un caract\`ere d'ordre fini tel que
$\varepsilon(-1)=(-1)^k$, de sorte que $w=(x \mapsto x^k\varepsilon(x)) \in \W(\Qpb)$.
On suppose que $k \neq 2$ si $\varepsilon=1$. Il existe 
(cf. \cite[thm. 4.7.1]{miyake}) une unique forme
modulaire classique, propre pour $\HH$, de $q$-d\'eveloppement 

$$E_w^{crit}:= q+\sum_{n\geq 2} a_n q^n, \, \, a_p=p^{k-1},\, \, a_l=\varepsilon(l)+l^{k-1} \, \, si \, \, \, l \neq
p.$$
Soit $F:=\Q_p(\varepsilon(\Z_p^*))$. La forme $E_w^{crit}$ d\'efinit donc un unique $F$-point de $\CC$
 que l'on note $x_w$, on a $\kappa(x_w)=w$. Il est clair que $\rho_{x_w}=\varepsilon
\oplus F(1-k)$. Un tel point de $\CC$ sera dit {\it Eisenstein critique}. \ps
	Soit $\zeta_p: \W\backslash\{1\} \longrightarrow \A^1$ la fonction z\^eta $p$-adique de
Kubota-Leopold. Si $w \neq 1$, Il existe une unique forme modulaire surconvergente, ordinaire
et propre pour $\HH$, de $q$-d\'eveloppement (cf. \cite[\S B1]{BMF}):

$$E_w^{\ord}:=\zeta_p(w)/2+q+\sum_{n\geq 2} a_n q^n, \, \, a_p=1,\, \,
a_l=1+w(l)l^{-1} \, \, si \, \, \, l \neq
p.$$

On pose de plus $E_1^{\ord}:=1$. Soit $F:=\Q_p(w(\Z_p^*))$. 
On notera $y_w$ le $F$-point de $\CC$ (en fait du lieu ordinaire $\CC^{\ord}$) correspondant. Un
point de la forme $y_w$ sera dit {\it Eisenstein ordinaire}. Les $y_w$ sont
en fait l'ensemble des points d'un ferm\'e Zariski $\CC^{\eis} \subset
\CC^{\ord}$, la {\it droite Eisenstein} (cf. \cite[\S 2.2]{eigen}), tel que $\kappa$ induit un isomorphisme $\CC^{\eis} \rightarrow
\W$.\ps
\medskip

{\bf Remarques:} \label{droiteeis} \ps 
i) Un point $y_w \in \CC^{\eis}$ est dans $\CC^0$ \ssi $\zeta_p(w)=0$. 
Ainsi, $\CC^{\eis}\cap \CC^0$ est non vide \ssi $p$ est un nombre
premier irr\'egulier. En g\'en\'eral, $\CC^{\eis}\cap \CC^0$ est fini
car $\zeta_p$ n'a qu'un nombre fini de z\'eros sur $\W$. \ps
ii) Rappelons (cf. \S \ref{composantes},\, \cite[\S 3.6]{eigen})) que $\CC^{\ord}=\CC^{\eis}\cup \CC^{0,\ord}$ est un ouvert ferm\'e de $\CC$,
et que $\CC^{0,\ord}$ est d'\'equidimension $1$. Comme $\CC^{\eis} \simeq \W$ est lisse, un point $x \in \CC^{\eis}$ est singulier vu comme
point de $\CC$ \ssi il est dans $\CC^0$, ou encore \ssi $\kappa$ est de
degr\'e $>1$ en $x$. V\'erifions que cela ne se produit
pas aux points classiques de $\CC^{\eis}$, i.e. aux $y_w$ tels que $w=(x
\mapsto x^k\varepsilon(x))$ avec $k\geq 1$ et $\varepsilon(-1)=(-1)^k$. Il suffit de v\'erifier que 
$\zeta_p(w)=L_p(1-k,\varepsilon)$ (cf. \cite[\S B1]{BMF} pour la notation) est non nul,
mais cela vient de ce que la fonction $L$ de
Dirichlet $L(s,\varepsilon)$ ne s'annule pas en $s=k$ si $k \geq 1$ et
$\varepsilon(-1)=(-1)^k$. 
\ps

\subsection{Points r\'eductibles de $\CC$.} Un point $x \in \CC(\Qpb)$ est dit
r\'eductible si $\rho_x$ l'est.\ps

\begin{prop} L'ensemble des points r\'eductibles de $\CC$ est exactement 
l'ensemble des points Eisenstein. 
\end{prop}

\pf Soit $x$ un point de $\CC(\Qpb)$ tel que $\rho_x=\chi_1+\chi_2$ est somme de deux
caract\`eres (automatiquement continus). D'apr\`es \cite[thm. 6.3]{kis}, pour $i=1$ ou
$2$, on a $D_{\mathrm{cris}}(\chi_i)^{\varphi=U_p(x)} \neq 0$. Supposons que $i=1$,
quitte \`a les renum\'eroter. Le caract\`ere $(\chi_1)_{|D}$ est donc 
cristallin de poids $k-1:=v(U_p(x)) \in \N$ (car $|U_p(x)| \leq 1$
pour tout $x$ dans $\CC$). D'autre part, $\chi_1$ est non ramifi\'e hors de $p$, c'est donc
$\Q_p(1-k)$. Comme les
poids de Hodge-Tate-Sen de $\rho_x$ sont $0$ et $d\kappa(x)-1$, il y a deux
possibilit\'es: \ps

\noindent - Si $v(U_p(x))>0$, alors $k\geq 2$. 
Il vient que $\varepsilon:=\chi_2$ est un caract\`ere de
poids $0$, donc d'ordre fini. Comme $E_2$ n'est pas surconvergente d'apr\`es
\cite{CGJ} (cf. aussi la discussion dans l'appendice de \cite{CheT}), 
$k=2 \Rightarrow \varepsilon \neq 1$. Ainsi, $x$ est de la forme $x_w$.
\ps

\noindent - Si $v(U_p(x))=0$, alors $\chi_1$ est le caract\`ere trivial et
donc $\chi_2=\det(\rho_x)$. Ainsi, $x$ est de la forme $y_w$. \epf

\medskip

La d\'emonstration ci-dessus montre de plus que

\begin{prop} \label{irrlocal} L'ensemble des $x \in \CC(\Qpb)$ tels que $v(U_p(x)) \neq 0$ et
$(\rho_x)_{|D}$ est r\'eductible est discret, compos\'e de $x$ tels que
$v(U_p(x))=d\kappa(x)-1$ est un entier strictement positif. 
\end{prop}

\section{Lissit\'e de $\CC$ aux points r\'eductibles non ordinaires}

\subsection{Rappels de cohomologie galoisienne.} Soient $k\geq 2$ un entier, 
$\varepsilon: \Z_p^* \rightarrow \Qpb^*$ un caract\`ere d'ordre fini tel que
$\varepsilon(-1)=(-1)^k$, et $F:=\Q_p(\varepsilon(\Z_p^*))$. On consid\`ere le caract\`ere de $G$  
$$\chi:=F(k-1)\otimes \varepsilon,$$
et on fait l'hypoth\`ese que $\chi \neq F(1)$. On rappelle le cas particulier suivant connu des conjectures de Bloch-Kato
(cf. \cite{BK}, \cite{FP}) pour
les fonctions $L$ de Dirichlet\footnote{Tous les groupes
de cohomologie galoisienne consid\'er\'es dans cette section sont
sous-entendu en cohomologie continue.}: 

$$\begin{array}{rclcc} \label{BKzeta} 
\dim_{F} H^1_f(\Q,\chi) & = & \ord_{s=2-k} L(s,\varepsilon^{-1}) & = & 1 \\
\dim_{F} H^1_f(\Q,\chi^{-1}) & = & \ord_{s=k} L(s,\varepsilon) & = & 0
\end{array}$$
Les \'egalit\'es de droite proviennent de l'\'equation fonctionnelle des
caract\`eres de Dirichlet, et de ce que $L(s,\varepsilon^{-1})$ n'a ni
z\'ero ni p\^ole en $s=n$ entier si $\varepsilon \neq 1$ et $n \geq 1$ ou si
$\varepsilon=1$ et $n\geq 2$. Les \'egalit\'es de gauche d\'ecoulent de
mani\`ere standard des travaux de Soulé (\cite{sou}). 
Notons que comme $\chi \neq \Q_p(1)$,
$H^1_f(\Q_p,\chi)=H^1(\Q_p,\chi)$ est de dimension $1$ et
$H^1_f(\Q_p,\chi^{-1})=0$. En particulier, $H^1_f(\Q,\chi)=H^1(G,\chi)$ et
$H^1_f(\Q,\chi^{-1})=\Ker(H^1(G,\chi^{-1}) \rightarrow H^1(\Q_p,\chi^{-1}))$.
Il est clair d'autre part que $H^1(G,\chi^{-1})$ est non nul (cela d\'ecoule
par exemple de la formule pour la caract\'eristique d'Euler
globale). En r\'ecapitulant, on obtient la:

\begin{prop} \label{sansf} \ps
i) Pour $H=G$ et $\Gal(\Qpb/\Q_p)$, on a 
$$\dim_F(\Ext^1_{cont, F[H]}(\varepsilon,
F(1-k)))=\dim_F(\Ext^1_{cont, F[H]}(F(1-k),\varepsilon))=1.$$ 
ii) L'application de restriction
$$\Ext^1_{cont,F[G]}(\varepsilon, F(1-k)) \longrightarrow
\Ext^1_{cont,F[\mathrm{Gal}(\Qpb/\Q_p)]}(\varepsilon,
F(1-k))$$ est un isomorphisme. \ps

\end{prop}

\subsection{} Dans ce qui suit, on se fixe un point Eisenstein critique
$x:=x_w$, $w: z \mapsto z^k\varepsilon(z)$. On a $x_w \in \CC(F)$ o\`u
$F=\Q_p(\varepsilon(\Z_p^*))$. On choisit un $F$-voisinage $\Omega$ de $x$ comme en
\S \ref{preparation} proposition \ref{voisinage}. 
Soit $A$ l'anneau local rigide de $\Omega$ en $x_w$ et nommons encore $T: G \rightarrow
A$ le pseudo-caract\`ere induit par $T$. L'anneau $A$ est une $F$-alg\`ebre
topologique (cf. \S \ref{exemple}) et $T: G \rightarrow A$ est continu par
d\'efinition de la topologie sur $A$ et \S \ref{Tcont}. \ps
Si $m$ d\'esigne l'id\'eal maximal de $A$, on a $A/m=F$ et par
construction $$T \bmod m= \varepsilon + F(1-k),$$ de sorte que $T$ est
r\'esiduellement somme de deux caract\`eres, distincts une fois restreints \`a $D$.
De plus, d'apr\`es la th\'eorie des pseudo-repr\'esentations de Wiles, il existe une repr\'esentation 
$\rho=(\rho_j): G \rightarrow \GL_2(K)$ de trace $T$, o\`u $K$ est l'anneau total de fractions de $A$
(cf. \S \ref{MWintro}) et les $\rho_j$ sont semi-simples. On pose $\chi_1:=F(1-k)$ et $\chi_2:=\varepsilon$. \ps
Fixons un $s \in D$ comme au \S \ref{MWintro} et choisissons une base
$s$-adapt\'ee de $K^2$. Cette base nous permet de d\'efinir $B$, $C$ (resp. $B_p$,
$C_p$) comme en \S \ref{MWintro} associ\'es \`a $\rho$ (resp. $\rho_{|D}$). On a
$B_p \subset B$ et $C_p \subset C$.  \ps

\begin{lemme} \label{defA} ${}^{}$ \ps

(a) Les $(\rho_j)_{|D}$ sont irr\'eductibles, \ps 

(b) $\rho$ est d\'efinie sur $A$, et quitte \`a changer de base adapt\'ee \`a $s$,
$B_p=B=A$. De plus, $C_p$ et $C$ sont libres
de rang $1$ sur $A$.
\end{lemme}

\pf V\'erifions le (a). Soient $A(\Omega_j)$ l'affinoide image de
$A(\Omega)$ dans $K_j$, $\Omega_j \subset \Omega$ le ferm\'e Zariski
contenant $x$ correspondant (il est d'\'equidimension $1$), et $T_j$ l'image de $T$ dans $A(\Omega_j)$. 
Si $(\rho_j)_{|D}$ est r\'eductible, les images de $B_p$
et $C_p$ dans $K_j$ sont nulles. Il vient que $(T_j)_{|D}$ est identiquement 
somme de deux caract\`eres, ce qui est absurde d'apr\`es la proposition
\ref{irrlocal}. \ps
D'apr\`es (a) et la proposition \ref{sansf}, le corollaire \ref{corMWtop}
(a) s'applique \`a $\rho$ et nous donne une base adapt\'ee \`a $s$ dans
laquelle $B=A$. En
particulier, $\rho$ est d\'efinie sur $A$. De m\^eme, il vient que $C$,
$B_p$ et $C_p$ sont libres de rang $1$ sur $A$. De plus, on a un diagramme commutatif 

$$\xymatrix{\Ext^1_{cont,F[G]}(\varepsilon, F(1-k)) \ar[r] & \Ext^1_{cont,F[\mathrm{\Gal}(\Qpb/\Q_p)]}(\varepsilon, F(1-k))
\\ \Hom_A(B,A/m) \ar@{->}^{\jj}[u] \ar[r]^{can} & \Hom_A(B_p,A/m)
\ar@{->}_{\jj}[u]}$$

\smallskip
\noindent La fl\`eche du haut est un isomorphisme entre $F$-espaces
vectoriels de dimension $1$ d'apr\`es la proposition \ref{sansf}, et les
fl\`eches verticales sont des isomorphismes d'apr\`es ce que l'on vient de
voir. Il vient que la fl\`eche du bas est un isomorphisme, ainsi donc que
l'inclusion $B_p \subset B$ d'apr\`es le lemme de Nakayama.\epf
\medskip 

\subsection{} On se place d\'efinitivement dans la base adapt\'ee donn\'ee par le lemme
ci-dessus, i.e. $B_p=B=A$. En particulier, $C=BC \subset m$ (resp.
$C_p=B_pC_p$) est un id\'eal de $A$: c'est l'id\'eal de r\'eductibilit\'e de
$T$ (resp. $T_{|D}$) en $x_w$. \ps

\begin{thm} \label{lissite} $A$ est de valuation discr\`ete et $C$ est l'id\'eal maximal de $A$. 
\end{thm}

\medskip

En particulier, $\CC$ est lisse en $x_w$. D'apr\`es le corollaire
\ref{corMWtop} (b), qui s'applique par le lemme \ref{defA} (a) et la proposition
\ref{sansf}, il suffit de montrer la seconde assertion. Un ingr\'edient crucial est la
cons\'equence suivante de \cite{kis}: \ps

\begin{lemme} \label{kis} Soit $J \subset m$ un id\'eal de codimension finie de $A$,
alors $D_{\mathrm{cris}}(\rho_{|D} \otimes A/J)^{\varphi=U_p}$ est libre de rang $1$
sur $A/J$. 
\end{lemme}

\pf Tout d'abord, notons que la trace de la repr\'esentation (\`a composantes
g\'en\'eriquement semi-simples) $\rho: G \rightarrow \GL_2(A)$ tombe dans l'anneau noeth\'erien
$A(\Omega)$. Par un argument standard d\'ej\`a donn\'e dans le lemme
\ref{tf}, $A(\Omega)[\rho(G)]$ est de type fini. Or 
$A$ est limite
inductive des $A(\Omega')$ o\`u $\Omega'$ parcourt les voisinages de $x$ de
dans $\CC$. Ainsi, quitte \`a r\'etr\'ecir $\Omega$ on peut
supposer que $\rho$ provient d'une repr\'esentation continue $$\rho^*: G \rightarrow
\GL_2(A(\Omega)),$$ par
$A(\Omega) \rightarrow A$. Quitte \`a restreindre encore $\Omega$, on peut
de plus supposer qu'il existe un id\'eal $J^*$ de $A(\Omega)$ tel
que $J^*A=J$ et l'application canonique $A(\Omega)/J^* \rightarrow A/J$ est
un isomorphisme. \ps
	Le polyn\^ome de Sen (cf. \cite{Se}) de $\rho^*$ est clairement de la forme
$T(T-d\kappa+1) \in A(\Omega)[T]$. On peut donc appliquer la proposition 5.4
de \cite{kis} \`a $X=\Omega$, $Y=U_p$ et $M=A(\Omega)^2$ muni de l'action de
$D$ via $\rho^*$. Par (2), $X_{fs}$ contient tous les points classiques de
$\Omega$, qui sont Zariski-denses dans $\Omega$. Comme c'est un ferm\'e
Zariski, $X_{fs}=\Omega$. On applique alors encore (2) et la remarque 5.5.
(1) \`a $f: \Max(A(\Omega)/J^*) \rightarrow \Omega$, 
qui est {\it $Y$-small} car de codimension finie, et se factorise par
$\Omega_{d\kappa-j}$ pour tout $j\leq 0$ car $d\kappa(x_w)-1=k-1 >0$. Cela
conclut. \epf
\ps
\ps

\subsection{} \label{findemo} Terminons la preuve. L'id\'eal $C$ est libre de rang $1$ dans l'anneau local
noeth\'erien $A$ qui est d'\'equidimension $1$, il est donc de codimension
finie d'apr\`es le Hauptidealsatz. Notons que $r:=\rho \otimes A/C$ est par
construction une extension 

$$ 0 \rightarrow (A/C).\psi_1 \rightarrow r \rightarrow (A/C).\psi_2
\rightarrow 0,$$
o\`u $\psi_1$ (resp. $\psi_2$) est la r\'eduction modulo $C$ de la fonction
$a$ (resp. $d$). Par construction, les $\psi_i: G \rightarrow (A/C)^*$ sont des caract\`eres
continus (cf. proposition \ref{topopo}), tels que 
$$\psi_1 \bmod m \equiv F(1-k), \, \, \, \, \, \psi_2 \bmod m \equiv
\varepsilon.$$

Consid\'erons la suite exacte de $A/C$-modules:

$$0 \rightarrow D_{\mathrm{cris}}((\psi_1)_{|D})^{\varphi=U_p} 
\rightarrow D_{\mathrm{cris}}(r_{|D})^{\varphi=U_p} \rightarrow
D_{\mathrm{cris}}((\psi_2)_{|D})^{\varphi=U_p} $$
Le lemme \ref{kis} implique que le terme central est libre de rang $1$.
Comme $\psi_2$, vu comme $F$-repr\'esentation, est une extension successive
de caract\`eres \'egaux \`a $\varepsilon$ et que
$$D_{\mathrm{cris}}(\varepsilon)^{\varphi=U_p(x)=p^{k-1}}=0,$$ le terme de droite est
nul. De plus, la dimension sur $F$ du terme de
gauche est $\leq \dim_F(A/C)$. Il vient donc que
$D_{\mathrm{cris}}((\psi_1)_{|D})^{\varphi=U_p}$ est libre de rang $1$ sur $A/C$. Cela implique que $(\psi_1)_{|D}(k-1)$ est 
cristallin de poids $0$, donc non ramifié. En particulier, le polyn\^ome de
Sen de $r_{|D}$ vaut $T(T-k+1)$, et donc $\kappa\equiv k \in A/C$, i.e.
$(\kappa-k) \subset C$. Comme le caractère $\psi_1(k-1)$ est d'autre part non ramifié hors de $p$,
il est identiquement trivial, puis  
$$\psi_1=(A/C)(1-k),\, \, \,  U_p \equiv p^{k-1} \in A/C,\, \, \, \mathrm{
et }\, \, \, \psi_2=\varepsilon.(A/C).$$
Le dernier point de la proposition \ref{voisinage} implique alors que $A/C=F$, i.e. $C$ est l'idéal maximal de $A$. \epf 

\smallskip

{\bf Remarques:} \ps
\noindent i) Une cons\'equence du lemme \ref{defA} (b) et de la premi\`ere
partie de la preuve du
lemme \ref{kis} est que le lieu non ordinaire de $\CC$ est admissiblement recouvert par des
ouverts affinoides
sur lesquels $T$ est la trace d'une vraie repr\'esentation de $G$ (c'est aussi une
cons\'equence simple du th\'eor\`eme \ref{lissite}, cf. la remarque
suivant \cite[thm. 5.1.2]{eigen}).  \ps

\noindent ii) Le th\'eor\`eme \ref{lissite} montre que
le diviseur de r\'eductibilit\'e de $T$ est r\'eduit. Est-il celui d'une
fonction globale sur $\CC$ ? \ps

\section{D\'etermination de $C_p$ et r\'egulateurs $p$-adiques}

\subsection{} Soit $x_w \in \CC(F)$ un point Eisenstein critique comme dans
la section pr\'ec\'edente. On reprend de plus les notations pr\'ec\'edentes pour $A$ et
$C_p$. Soient $\chi=F(k-1)\otimes \varepsilon$ le caract\`ere de $G$
correspondant \`a $w \in \W(F)$, et $w^*:=(z \mapsto z^{2-k}\varepsilon^{-1}(z))$ le poids du
point Eisenstein ordinaire $y_{w^*}$ de $\CC(F)$ jumeau \`a $x_w$. \ps

\begin{thm} \label{calculCp} $C_p=(\kappa-k)$ et $\kappa$ a m\^eme degr\'e en $x_w$
qu'en $y_{w^*}$. \end{thm}

\pf La d\'emonstration du \S \ref{findemo} ci-dessus appliqu\'ee \`a $C_p$
plut\^ot qu'\`a $C$ montre encore que $(\kappa-k) \subset C_p$. Nous avons
donc que $(\kappa-k) \subset C_p=B_pC_p$, ce dernier \'etant
l'id\'eal de $A$ de r\'eductibilit\'e de $T_{|D}$. D'apr\`es la remarque
d\'ebutant le \S \ref{defidred}, il suffit donc de montrer que $T: D \rightarrow A/(\kappa-k)=\HH(x_w)$ est somme de deux
caract\`eres. Le cas limite du crit\`ere de classicit\'e de Coleman \cite[cor. 7.2.2]{CO},
\cite{CO2}, montre que l'op\'erateur $\theta^{k-1}$ induit un isomorphisme $\HH\otimes_{\Z}F$-\'equivariant:
\begin{equation} \label{etoile} M(y_{w^*})\otimes_F \nu^{k-1}
\overset{\sim}{\rightarrow} M(x_w),
\end{equation}
o\`u $\nu: \HH \rightarrow F$ est le morphisme d'anneaux d\'efini par $T_l \mapsto l, \, \,  U_p \mapsto p$.
En particulier, vue la proposition \ref{voisinage} b), cela prouve la seconde assertion du th\'eor\`eme et montre
que l'on dispose d'un isomorphisme de
$F$-alg\`ebres locales $$\HH(y_{w*}) \overset{\sim}{\rightarrow}
\HH(x_w),$$ qui est un morphisme de $\HH$-alg\`ebres s'il on tord
l'application naturelle $\HH \rightarrow \HH(y_{w*})$ par $\nu^{k-1}$. 
D'apr\`es le th\'eor\`eme de Cebotarev, on en d\'eduit que via l'identification
ci-dessus, les pseudo-caract\`eres d\'eduits de $T$ par \'evaluations, $T: G \rightarrow \HH(x_w)$ et
$T(1-k): G \rightarrow \HH(y_{w^*})$, sont \'egaux. Mais sur le lieu ordinaire
$\CC^{\ord} \subset \CC$ de $\CC$, qui est un ouvert admissible de $\CC$ contenant
$y_{w^*}$, $T_{|D}$ est r\'eductible. En effet, il est m\^eme somme du caract\`ere non ramifi\'e
$\chi: D \rightarrow A(\CC)^*$ envoyant un Frobenius g\'eom\'etrique sur $U_p$
et de $\chi^{-1}\det(T)$. Cela conclut.
\epf

\ps

En particulier, nous avons montr\'e le:

\begin{cor} \label{reploc} La restriction \`a $D$ de la repr\'esentation 
$\rho: G \rightarrow \GL_2(\HH(x_w))$ est une
extension de  $\HH(x_w)(\varepsilon \otimes u^{-1})$ par $\HH(x_w)(1-k)\otimes u$, o\`u
$u$ est le caract\`ere non ramifi\'e $D \rightarrow \HH(x_w)^*$ envoyant le
Frobenius g\'eom\'etrique sur $U_p/p^{k-1}$. 
\end{cor}

\subsection{} Il ne semble cependant pas possible d'en d\'eduire que $\rho$ est constante,
i.e. que $C_p$ est l'id\'eal maximal de $A$. Cela vient de ce que l'on ne sait pas si le morphisme
r\'egulateur $H^1(G,\chi) \rightarrow H^1(\Q_p,\chi)$ est trivial ou non. 
Concernant ce morphisme on a en fait le :
\ps

\begin{thm} Les propri\'et\'es suivantes sont \'equivalentes: \ps

i) $\kappa$ est \'etale en $x_w$, \ps
i') $\dim_F M(x_w)=1$,\ps
ii) L'application naturelle $H^1(G,\chi) \rightarrow H^1(\Q_p,\chi)$ est un isomorphisme, \ps
iii) $\zeta_p(w^*)\neq 0$. \ps

\end{thm}

\pf On a d\'ej\`a vu que i) est \'equivalent \`a i') (cf. proposition
\ref{voisinage} b)). De plus, on montre comme dans le lemme \ref{defA} que
ii) est \'equivalent \`a ce que $C_p=C$. Mais ceci est \'equivalent 
\`a ce que $(\kappa-k)$ soit l'id\'eal maximal de $A$ d'apr\`es les th\'eor\`emes
\ref{lissite} et \ref{calculCp}. Cela montre
l'\'equivalence de i) et ii). L'\'equivalence entre ii) et iii) est bien
connue, nous allons la red\'emontrer ici en v\'erifiant celle de i) et iii).
D'apr\`es (\ref{etoile}), $\kappa$ est de degr\'e $1$ (i.e. \'etale) en $x_{w^*}$ \ssi
$\kappa$ est de degr\'e $1$ en $y_{w^*} \in \CC^{\eis}$. Mais d'apr\`es les
remarques \ref{droiteeis} i) et ii), ceci se produit \ssi $\zeta_p(w^*)\neq 0$.
\epf


\medskip

Il est commun\'ement conjectur\'e que les propri\'et\'es ii) et iii) \'equivalentes
ci-dessus sont satisfaites, ainsi donc que i), i'). Si $p$ est un nombre
premier r\'egulier, elles sont toujours satisfaites. \ps

\end{document}